[1994/12/01]
\documentclass[draft]{article}
\title{The Extension of Knot Groups to Tangles}
\author{John Armstrong}
\date{September 28, 2005}

\usepackage{amsmath,amsthm,amssymb}
\usepackage{amscd}
\usepackage{graphicx}
\usepackage[all]{xy}
\xyoption{knot}
\CompileMatrices

\input{epsf}

\newtheorem{Theorem}{Theorem}
\newtheorem{Corollary}{Corollary}
\newtheorem{Conjecture}{Conjecture}

\theoremstyle{definition}
\newtheorem{Definition}{Definition}

\begin{document}
\maketitle

\begin{abstract}
  The extension of the knot group $\pi_1(S^3\setminus K)$ to the category of tangles is introduced via a new category-theoretic construction.  Through this presentation, a new avenue of proof for results about knot groups is opened.
\end{abstract}

\renewcommand{\sectionmark}[1]{}

\section{Introduction}

One of the better-known invariants of knots and links is the ``knot group'' $\Gamma(K)=\pi_1(S^3\setminus K)$, defined as the fundamental group of the complement of a knot or link in $S^3$~\cite{MR0176462}.  One would like to extend this invariant to tangles.  That is: to define a functor $\Gamma$ from $\mathcal{T}ang$ -- the category of tangles -- to an appropriate target category so that its restriction to $\hom_{\mathcal{T}ang}(0,0)$ -- knots and links -- is essentially the same as the knot group function.  Since $\mathcal{T}ang$ has a simple presentation generators and relations as a tensor category, this will allow a new computation of the knot group of a given knot or link, and more importantly a new avenue of proof for results on knot groups.

The na\"ive choice of $\mathbf{Grp}$ for the target category is obviously incorrect.  Groups are objects in $\mathbf{Grp}$, while tangles are morphisms in $\mathcal{T}ang$.  As a first step, then, one must construct a category in which (essentially) groups are the morphisms.  

\section{The $\Lambda$ construction}

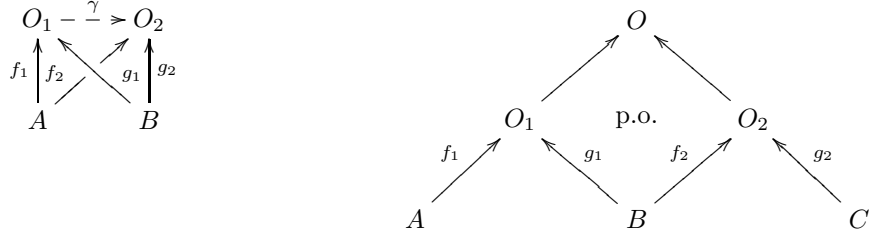
\begin{figure}
  \centerline{
    \xymatrix{
      O_1 \ar@{-->}[r]^{\gamma} & O_2\\
      A \ar[u]^{f_1} \ar[ur]|{\hole}^(.3){f_2} & B \ar[ul]_(.3){g_1} \ar[u]_{g_2}
    }
    \qquad\qquad\qquad\qquad
    \xymatrix{
      && O &&\\
      & O_1 \ar[ur] & \text{p.o.} & O_2 \ar[ul] &\\
      A \ar[ur]^{f_1} && B \ar[ul]_{g_1} \ar[ur]^{f_2} && C \ar[ul]_{g_2}
    }
  }
  \caption{Isomorphism and composition of diagrams}
  \label{fig:isomcomp}
\end{figure}

Let $\mathcal{C}$ be a category with pushouts.  Construct a new category $\Lambda\mathcal{C}$ with the same objects as $\mathcal{C}$, but with new morphism sets
\begin{equation*}
\hom_{\Lambda\mathcal{C}}(A,B) = \{\text{isomorphism classes of diagrams } A \xrightarrow{f} O \xleftarrow{g} B\}
\end{equation*}
where two diagrams $A \xrightarrow{f_1} O_1 \xleftarrow{g_1} B$ and $A \xrightarrow{f_2} O_2 \xleftarrow{g_2} B$ are considered isomorphic if there is an isomorphism $\gamma$ from $O_1$ to $O_2$ such that $f_2 = f_1\circ\gamma$ and $g_1 = \gamma\circ g_2$, as in Figure \ref{fig:isomcomp}.

Given a morphism from $A$ to $B$ and one from $B$ to $C$, define their composition by splicing diagrams at $B$ and pushing out $g_1$ and $f_2$ as in Figure \ref{fig:isomcomp}.  Pushouts are only defined uniquely up to isomorphism, which is why morphisms in $\Lambda\mathcal{C}$ must be defined as isomorphism classes of diagrams.

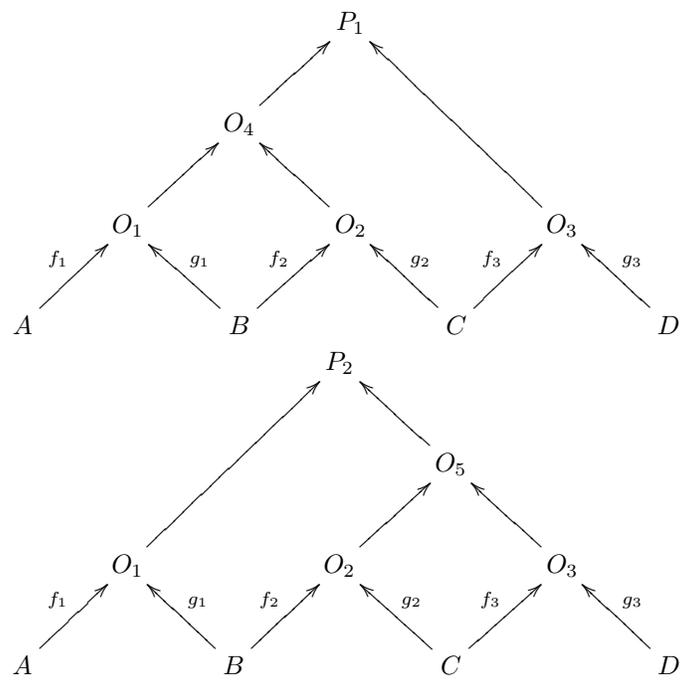
\begin{figure}
  \centerline{
    \xymatrix{
      &&& P_1 &&&\\
      && O_4 \ar[ur] &&&&\\
      & O_1 \ar[ur] && O_2 \ar[ul] && O_3 \ar[uull] &\\
      A \ar[ur]^{f_1} && B \ar[ul]_{g_1} \ar[ur]^{f_2} && C \ar[ul]_{g_2} \ar[ur]^{f_3} && D \ar[ul]_{g_3}
    }
  }
  \centerline{
    \xymatrix{
      &&& P_2 &&&\\
      &&&& O_5 \ar[ul] &&\\
      & O_1 \ar[uurr] && O_2 \ar[ur] && O_3 \ar[ul] &\\
      A \ar[ur]^{f_1} && B \ar[ul]_{g_1} \ar[ur]^{f_2} && C \ar[ul]_{g_2} \ar[ur]^{f_3} && D \ar[ul]_{g_3}
    }
  }
  \caption{Associativity of Composition}
  \label{fig:assoc}
\end{figure}

Considering three morphisms and both orders of composition, it is easily verified that both compositions $P_1$ and $P_2$ (see Figure \ref{fig:assoc}) are limits of the diagram composed of the original morphisms, and so are isomorphic. 
The identity diagram for an object $A$ is $A \xrightarrow{id_A} A \xleftarrow{id_A} A$.

\begin{Theorem}
  \label{thm:subcat}
  There exists a faithful functor $\lambda$ from $\mathcal{C}$ to $\Lambda\mathcal{C}$
\end{Theorem}
\begin{proof}
  Set $\lambda(A) = A$ and $\lambda(A\xrightarrow{f}B) = {A \xrightarrow{f} B \xleftarrow{id_B} B}$.  This is clearly a functor, so it remains to show faithfulness.  But if two diagrams from $A$ to $B$ in the image of $\lambda$ are isomorphic, the isomorphism on the middle object must be the identity, and so the arrows from $A$ must be identical.
\end{proof}

Now let $\otimes$ be a tensor product (associative bifunctor with identity object $\mathbf{1}$ up to natural isomorphism) on $\mathcal{C}$ which commutes with pushouts.  This extends naturally to a tensor product on $\Lambda\mathcal{C}$.  Indeed, $A_1\otimes A_2 \xrightarrow{f_1\otimes f_2} O_1\otimes O_2 \xleftarrow{g_1\otimes g_2} B_1\otimes B_2$ is a morphism from $A_1\otimes A_2$ to $B_1\otimes B_2$, and the pushout of the tensor product of two arrows is the tensor product of their pushouts by assumption, proving functoriality.  The identity diagram on $\mathbf{1}$ acts as a tensor identity in $\Lambda\mathcal{C}$, and associativity is easily checked.

Conversely, given any tensor product structure on $\Lambda\mathcal{C}$, we can recover one on $\mathcal{C}$ by restricting to the image of $\lambda$.  It evidently preserves pushouts.  This proves

\begin{Theorem}
  \label{thm:tensstruct}
  Tensor product structures on $\Lambda\mathcal{C}$ correspond bijectively with tensor product structures on $\mathcal{C}$ that preserve pushouts.
\end{Theorem}

\section{The functor $\Gamma$}

The category $\mathbf{Grp}$ has pushouts, as well as an associative (up to natural isomorphism), unital (again, up to isomorphism) bifunctor preserving them: the coproduct, better known as the free product.  This will induce a tensor structure on $\Lambda\mathbf{Grp}$, and this structure will be assumed in the sequel.  The tensor identity is the single element group $\mathbf{1}$.

$\mathcal{T}ang$ has as objects the natural numbers with addition acting as a monoidal product.  A tangle from $m$ to $n$ is an isotopy class of collections of $\frac{m+n}{2}$ arcs and any number of circles embedded into the slab $\mathbb{R}^2\times[0,1]$.  The endpoints of the arcs lie at the points $(1,0),\dots,(m,0)$ on the lower face of the slab and $(1,0),\dots,(n,0)$ on the upper face.  Given a tangle from $m$ to $n$ and one from $n$ to $p$, they are composed by stacking their slabs atop each other to line up the $n$ free ends from each tangle and then reparametrizing the height of the slab to $[0,1]$.  Finally, given tangles $T_1\in\hom_{\mathcal{T}ang}(m_1,n_1)$ and $T_2\in\hom_{\mathcal{T}ang}(m_2,n_2)$, the tensor product $T_1\otimes T_2$ is given  by sliding the lower free ends of $T_2$ to $(m_1+1,0),\dots,(m_1+m_2,0)$, the upper free ends to $(n_1+1,0),\dots,(n_1+n_2,0)$, and embedding $T_1$ and $T_2$ into the same slab such that there exists a strip in the slab separating the image of $T_1$ from that of $T_2$.

To preserve the tensor product structure, $\Gamma(0) = \mathbf{1}$ and $\Gamma(n) = \Gamma(1)^{*n}$.  Note already what this implies for knots and links: $\Gamma$ will send a knot to a diagram of the form $\mathbf{1} \xrightarrow{f} G \xleftarrow{g} \mathbf{1}$.  But as $\mathbf{1}$ is initial the morphisms $f$ and $g$ are uniquely determined.  The only information left to determine is the group $G$ itself, and this is where the classical knot group $\Gamma(K)$ fits into the new picture.

\begin{Definition}
  \label{def:Gammafunctor}
  For any tangle $T\in\hom_{\mathcal{T}ang}(m,n)$
  \begin{multline*}
    \Gamma(T) =\\
    \pi_1(\mathbb{R}^2\setminus\{(1,0),...,(m,0)\}) \xrightarrow{f} \pi_1(\mathbb{R}^2\times [0,1]\setminus T) \xleftarrow{g} \pi_1(\mathbb{R}^2\setminus\{(1,0),...,(n,0)\})
  \end{multline*}
  where $f$ and $g$ are the homomorphisms induced by the inclusions of the boundary planes into the slab $\mathbb{R}^2\times [0,1]$.
\end{Definition}

\begin{Theorem}
 \label{thm:functoriality}
 $\Gamma$ is a tensor functor from $\mathcal{T}ang$ to $\Lambda\mathbf{Grp}$
\end{Theorem}
\begin{proof}
  The identity tangle on $n$ is the collection of $n$ line segments connecting $(i,0,0)$ to $(i,0,1)$ in the slab.  The fundamental group of its complement is $F_n$, and the inclusions of the boundary planes induce the identity morphism on fundamental groups, as desired.
  
  Given $T_1\in\hom_{\mathcal{T}ang}(m,n)$ and $T_2\in\hom_{\mathcal{T}ang}(n,p)$, the plane between the slabs in the composition may be thickened slightly to form a collar also having fundamental group $F_n$.  The fundamental group of the complement of $T_1\circ T_2$ may be calculated by the Seifert-van Kampen theorem as the product of the fundamental groups of the complements of $T_1$ and $T_2$ amalgamated over the fundamental group of this collar.  In categorical terms, this is the pushout of
  \begin{equation*}
    \pi_1(\mathbb{R}^2\times [0,1]\setminus T_1)\leftarrow F_n\rightarrow\pi_1(\mathbb{R}^2\times [0,1]\setminus T_2)
  \end{equation*}
  and the inclusions of the lower boundary of $T_1$ and the upper boundary of $T_2$ into the slab induce the compositions of their original inclusion morphisms with the pushout morphisms.  Thus $\Gamma(T_1\circ T_2) = \Gamma(T_1)\circ\Gamma(T_2)$, as desired (Figure \ref{fig:SvK}).
  
  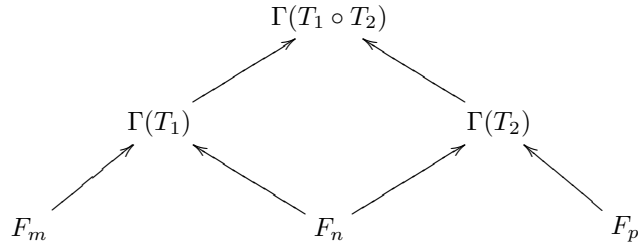
\begin{figure}[b]
  \centerline{
    \xymatrix{
      && \Gamma(T_1\circ T_2) &&\\
      & \Gamma(T_1) \ar[ur] && \Gamma(T_2) \ar[ul] &\\
      F_m \ar[ur] && F_n \ar[ul] \ar[ur] && F_p \ar[ul]
    }
  }
  \caption{Composing two Fundamental Groups}
  \label{fig:SvK}
\end{figure}
  
  On objects, $\Gamma(m+n) = F_{m+n} = F_m*F_n = \Gamma(m)*\Gamma(n)$, so $\Gamma$ indeed preserves the tensor product on objects.
  
  Finally, in the tensor product of two tangles, the separating strip may be thickened slightly.  The fundamental group of the complement of $T_1\otimes T_2$ may again be calculated by the Seifert-van Kampen theorem as the product of the complements of $T_1$ and $T_2$ amalgamated over the fundamental group of this thickened strip, which of course is trivial.  Thus $\Gamma(T_1\otimes T_2) = \Gamma(T_1)*\Gamma(T_2)$, as desired.
\end{proof}

Since the fundamental group of a knot or link complement in the slab is isomorphic to the fundamental group of its complement in $S^3$, we immediately have

\begin{Corollary}
  \label{cor:restknotlink}
  For a knot or link $K\in\hom_{\mathcal{T}ang}(0,0)$,
  \begin{equation*}
    \Gamma(K) = \mathbf{1}\xrightarrow{1}G\xleftarrow{1}\mathbf{1}
  \end{equation*}
  where $G$ is the knot group of $K$.
\end{Corollary}

\begin{figure}[h]
  \centerline{
    \mbox{
      \begin{xy}
        0;/r1cm/:,\vcrossneg
      \end{xy}
    }
    \qquad\qquad
    \mbox{
      \begin{xy}
        0;/r1cm/:,\vcross
      \end{xy}
    }
    \qquad\qquad
    \mbox{
      \begin{xy}
        0;/r1cm/:,\vcap-
      \end{xy}
    }
    \qquad\qquad
    \mbox{
      \begin{xy}
        0;/r1cm/:-(0,1),\vcap
      \end{xy}
    }
  }
  \caption{Generating tangles}
  \label{fig:tanggen}
\end{figure}
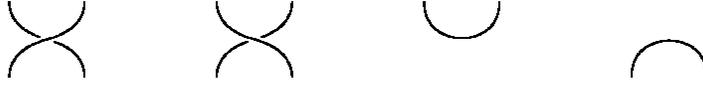

The category of tangles is isomorphic to the category of tangle diagrams, which has a simple presentation by generators and relations.  Its class of objects is the set of natural numbers with addition as a tensor product.  The morphisms are generated by the right-handed and left-handed crossings, cup, and cap (see Figure \ref{fig:tanggen}), which will be denoted $X^+$, $X^-$, $\cup$, and $\cap$, respectively.

These are subject to the relations

\begin{align*}
  T0:& (\cup\otimes id_1)\circ(id_1\otimes\cap) = id_1 = (id_1\otimes\cup)\circ(\cap\otimes id_1)\\
  T0':& (id_1\otimes\cup)\circ(X^\pm\otimes id_1) = (\cup\otimes id_1)\circ(id_1\otimes X^\mp)\\
  T0'':& (X^\pm\otimes id_1)\circ(id_1\otimes\cap) = (id_1\otimes X^\mp)\circ(\cap\otimes id_1)\\
  T1:& (id_1\otimes\cap)\circ(X^\pm\otimes id_1)\circ(id_1\otimes\cup)\\
  T2:& X^\pm\circ X^\mp = id_2\\
  T3:& (X^+\otimes id_1)\circ(id_1\otimes X^+)\circ(X^+\otimes id_1) =\\
  &\qquad\qquad\qquad\qquad= (id_1\otimes X^+)\circ(X^+\otimes id_1)\circ(id_1\otimes X^+)\\
\end{align*}

This is essentially the same construction given by Freyd and Yetter~\cite{MR1020583}.  The relations $T1$, $T2$, and $T3$ correspond to the Reidemeister moves of types I, II, and III, respectively, while relations $T0$, $T0'$, and $T0''$ handle other plane isotopies.

\begin{figure}[h]
  \epsfxsize=\textwidth
  \epsfbox{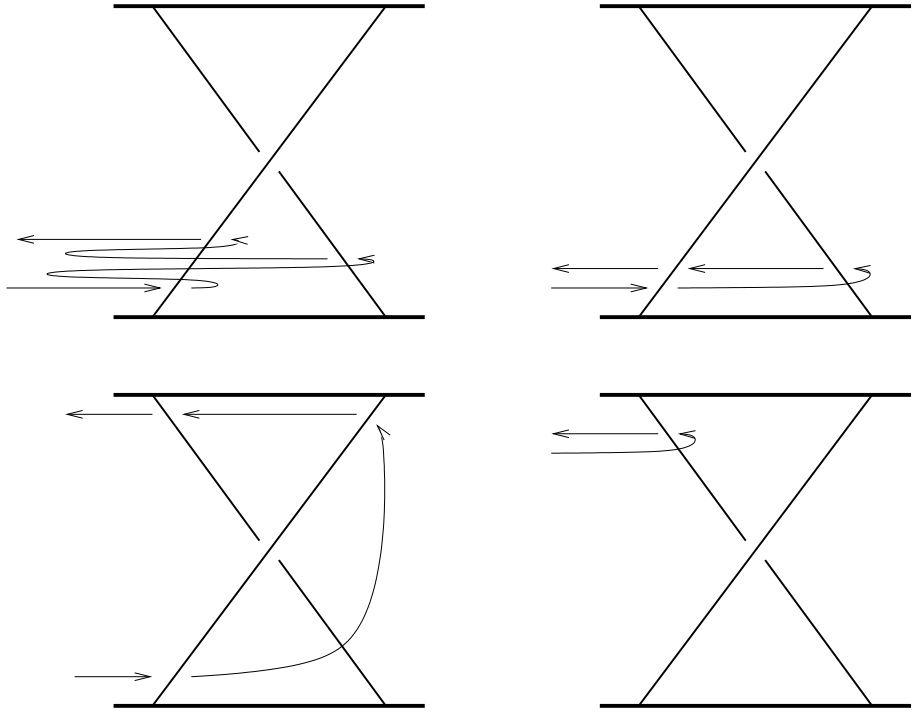}
  \caption{$\Gamma(X^+)$.  $c_1 = a_1^{-1}a_2a_1$}
  \label{fig:gammax}
\end{figure}

\begin{Theorem}
  \label{thm:generators}
  \begin{align*}
    \Gamma(\cup) &= \mathbf{1} \xrightarrow{1} <b> \xleftarrow{a_1\mapsto b, a_2\mapsto b^{-1}} <a_1,a_2>\\
    \Gamma(\cap) &= <a_1,a_2> \xrightarrow{a_1\mapsto b, a_2\mapsto b^{-1}} <b> \xleftarrow{1} \mathbf{1}\\
    \Gamma(X^+) &= <a_1,a_2> \xrightarrow{a_1\mapsto b_1, a_2\mapsto b_2} <b_1,b_2> \xleftarrow{c_1\mapsto b_1^{-1}b_2b_1,c_2\mapsto b_1} <c_1,c_2>\\
    \Gamma(X^-) &= <a_1,a_2> \xrightarrow{a_1\mapsto b_1, a_2\mapsto b_2} <b_1,b_2> \xleftarrow{c_1\mapsto b_2,c_2\mapsto b_2b_1b_2^{-1}} <c_1,c_2>\\
  \end{align*}
\end{Theorem}
\begin{proof}
  In each case, the generators of $\pi_1(\mathbb{R}^2\setminus\{(1,0),...,(m,0)\})$ are chosen with basepoint $(0,0)$, passing along the front of the first $i$ punctures (viewed from the $-y$ axis), around the back of the $i$th puncture, and back along the front of the first $i-1$ punctures.
  
  For the cup and cap, $\mathbb{R}^2\times [0,1]\setminus T$ is a slab with a single tube removed, so its fundamental group is $F_1$.  The two generators of the doubly-pierced boundary plane wrap in opposite directions around the tube, so one goes to the generator of $\pi_1(\mathbb{R}^2\times [0,1]\setminus T)$ and the other to its inverse.
  
  $\mathbb{R}^2\times [0,1]\setminus X^\pm$ is a slab with two tubes removed, so its fundamental group is $F_2$.  Choose the two generators to be those from the lower boundary plane.  The inclusion of the generator from the upper boundary plane encircling the overcrossing arc is then identical to the inclusion of the generator from its other end, while the inclusion of the generator encircling the undercrossing arc is the conjugate of the inclusion of the generator from its other end (Figure \ref{fig:gammax}).
\end{proof}


\section{Inductive Proofs}

Since $\mathcal{T}ang$ is a tensor category with a known set of generators, any tangle can be constructed starting with generators and applying compositions and tensor products.  Accordingly, statements about tangles can be proven by an inductive technique.  While this is just a special case of induction on a construction, it is useful to state it explicitly here.

\begin{Theorem}
  \label{thm:induction}
  Given a predicate $P$ on tangles, if the following are true
  \begin{itemize}
    \item $P$ is true for $id_1\in\hom_{\mathcal{T}ang}(1,1)$
    \item $P$ is true for $\cup$
    \item $P$ is true for $\cap$
    \item $P$ is true for $X^+$
    \item $P$ is true for $X^+$
    \item If $P$ is true for tangles $T_1\in\hom_{\mathcal{T}ang}(m,n)$ and $T_2\in\hom_{\mathcal{T}ang}(n,p)$ then it is true for $T_1\circ T_2$
    \item If $P$ is true for any two tangles $T_1$ and $T_2$, then it is true for $T_1\otimes T_2$
  \end{itemize}
  then $P$ is true for all tangles.  In particular, it is true for all knots and links.
\end{Theorem}

As an example of the use of this technique, it is possible to reconstruct a well-known theorem on knot groups.

\begin{Theorem}
  Let $T\in\hom_{\mathcal{T}ang}(m,n)$ be a tangle with $l$ loops and let $\Gamma(T) = F_m\rightarrow G\leftarrow F_n$.  Then the abelianization of $G$ is free on $\frac{m+n}{2}+l$ generators.
\end{Theorem}
\begin{proof}
  Each generating tangle may be easily checked to satisfy the conditions of the theorem.
  
  Let $T_1\in\hom_{\mathcal{T}ang}(m_1,n_1)$ and $T_2\in\hom_{\mathcal{T}ang}(m_2,n_2)$ with $l_1$ and $l_2$ loops, respectively.  Assume $\Gamma(T_1) = F_{m_1}\rightarrow G_1\leftarrow F_{n_1}$ where ${\rm Ab}(G_1)$ is free on $\frac{m_1+n_1}{2}+l_1$ generators, and similarly for $\Gamma(T_2)$.  Now
  \begin{equation*}
    \Gamma(T_1\otimes T_2) = F_{m_1+m_2}\rightarrow G_1*G_2\leftarrow F_{n_1+n_2}
  \end{equation*}
  
  Since ${\rm Ab}(G_1*G_2) = {\rm Ab}(G_1)\oplus{\rm Ab}(G_2)$, it is free on $\frac{(m_1+m_2)+(n_1+n_2)}{2}+l_1+l_2$ generators.  Since $T_1\otimes T_2\in\hom_{\mathcal{T}ang}(m_1+m_2,n_1+n_2)$ with $l_1+l_2$ loops, this satisfies the condition of the theorem.
  
  Let $T_1\in\hom_{\mathcal{T}ang}(m,n)$ and $T_2\in\hom_{\mathcal{T}ang}(n,p)$ with $l_1$ and $l_2$ loops, respectively.  Assume $\Gamma(T_1) = F_m\rightarrow G_1\leftarrow F_n$ where ${\rm Ab}(G_1)$ is free on $\frac{m+n}{2}+l_1$ generators, and similarly for $\Gamma(T_2)$.  Now
  \begin{equation*}
    \Gamma(T_1\circ T_2) = F_m\rightarrow G_1\times_{F_n}G_2\leftarrow F_p
  \end{equation*}
  
  Since the abelianization functor preserves pushouts, the abelianization of $G_1\times_{F_n}G_2$ is ${\rm Ab}(G_1)\oplus {\rm Ab}(G_2)$ modulo $n$ additional relations.  Each such relation comes from splicing a free end of $T_1$ to one of $T_2$, identifying a generator of ${\rm Ab}(G_1)$ with one of ${\rm Ab}(G_2)$.  However, when a given splicing closes a new loop the other ends have already been identified, so the new relation is redundant.  Thus if $l$ new loops are introduced only $n-l$ generators are killed.  Counting up, ${\rm Ab}(G_1\times_{F_n}G_2)$ is thus free on
  \begin{equation*}
    \frac{m+n}{2}+l_1+\frac{n_p}{2}+l_2-(n-l) = \frac{m+p}{2}+l_1+l_2+l
  \end{equation*}
  generators.  Since $T_1\circ T_2\in\hom_{\mathcal{T}ang}(m,p)$ and has $l_1+l_2+1$ loops -- $l_1+l_2$from $T_1$ and $T_2$ and $l$ created in the composition -- this satisfies the condition of the theorem.
  
  Thus, since the theorem holds for the generating tangles and for the two rules of composition, it holds for all tangles.
\end{proof}

Now, if $T$ is a link there is the special case

\begin{Corollary}
  The abelianization of the group of an $n$-component link $L$ is free on $n$ generators.
\end{Corollary}

This sort of technique may provide inductive algorithms for calculation not only of $\Gamma$, but of other generalized invariants related to the knot complement.  For instance, if the A-Polynomial~\cite{MR1288467} can be suitably generalized, it would immediately come with an inductive algorithm for its calculation.

\section{Skein Theory}

The notion of a skein theory for an invariant of links was introduced by Conway~\cite{MR0258014}, but the term has as yet no good definition beyond ``I know it when I see it''~\cite{SCOTUS378US184}.  Roughly it is used to localize information about the invariant; changing this part of a link this way changes the value of the invariant that way.  As such, most topologists would agree that $\Gamma$ is inherently a global notion and cannot be so localized.  Indeed the definition seems inherently global, but the Gauss integral also seems inherently global and yet there is a well-known skein theory defining the linking number~\cite{MR0515288}.

To probe this further, what is needed is a

\begin{Definition}
  A skein theory for for an invariant $F\colon\mathcal{T}ang\to\mathcal{C}$ is a finitely-generated tensor ideal $\mathcal{I}$ of $\mathcal{T}ang$ such that $F$ factors through the canonical projection to $\mathcal{T}ang/\mathcal{I}$.
\end{Definition}

Evidently a skein theory is not uniquely defined for an invariant of tangles since the projection onto $\mathcal{T}ang/\mathcal{I}$ itself factors through the quotient of $\mathcal{T}ang$ by any smaller ideal.  If $\tilde{F}\colon\mathcal{T}ang/\mathcal{I}\to\mathcal{C}$ is faithful, however, $\mathcal{I}$ is evidently maximal.  When considering functors generalized from link invariants it is convenient to require only that the restriction of $\tilde{F}$ to the images of knots and links in $\mathcal{T}ang/\mathcal{I}$ be injective.

\begin{Conjecture}
  There is no skein theory for $\Gamma$.  That is: there is no finitely-generated tensor ideal $\mathcal{I}$ of $\mathcal{T}ang$ so that $\Gamma$ factors through $\mathcal{T}ang/\mathcal{I}$, and also $\tilde{\Gamma}\colon\mathcal{T}ang/\mathcal{I}\to\Lambda\mathbf{Grp}$ is injective.
\end{Conjecture}

\begin{figure}
  \epsfxsize=\textwidth
  \epsfbox{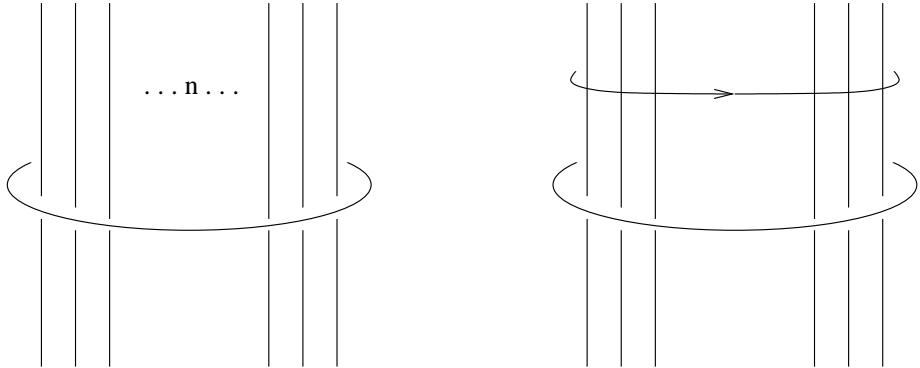}
  \caption{An infinite family of pairs of tangles not distinguished by $\Gamma$.}
  \label{fig:kernGamma}
\end{figure}

Unfortunately there is as yet no proof of this assertion, though it seems to be widely believed.  $\Gamma$ is known to be a perfect invariant for knots, but not for links~\cite{MR0515288}.  For example, the two tangles in Figure \ref{fig:kernGamma} can be calculated to have the same value of $\Gamma$.  Is there a finite collection of relations which implies this equality for all values of $n$ and which does not imply unwanted relations?  The intuition is against it, but it is admittedly unknown.  At least now the question can be stated sensibly and a proper attack may be made.


\pagebreak
\bibliographystyle{plain}
\bibliography{biblio}

\begin{thebibliography}{1}

\bibitem{MR0258014}
J.~H. Conway.
\newblock An enumeration of knots and links, and some of their algebraic
  properties.
\newblock In {\em Computational Problems in Abstract Algebra (Proc. Conf.,
  Oxford, 1967)}, pages 329--358. Pergamon, Oxford, 1970.

\bibitem{MR1288467}
D.~Cooper, M.~Culler, H.~Gillet, D.~D. Long, and P.~B. Shalen.
\newblock Plane curves associated to character varieties of {$3$}-manifolds.
\newblock {\em Invent. Math.}, 118(1):47--84, 1994.

\bibitem{MR1020583}
Peter~J. Freyd and David~N. Yetter.
\newblock Braided compact closed categories with applications to
  low-dimensional topology.
\newblock {\em Adv. Math.}, 77(2):156--182, 1989.

\bibitem{MR0176462}
L.~P. Neuwirth.
\newblock {\em Knot groups}.
\newblock Annals of Mathematics Studies, No. 56. Princeton University Press,
  Princeton, N.J., 1965.

\bibitem{SCOTUS378US184}
Supreme~Court of~the United~States.
\newblock Jacobellis v. {O}hio.
\newblock 378 U.S. 184 (1964).

\bibitem{MR0515288}
Dale Rolfsen.
\newblock {\em Knots and links}.
\newblock Publish or Perish Inc., Berkeley, Calif., 1976.
\newblock Mathematics Lecture Series, No. 7.

\end{thebibliography}
\end{document}